\documentclass[a4paper,12pt]{amsart}
\usepackage{amsmath,amssymb,latexsym,amsfonts,amscd}

\title{The canonical volume of 3-folds of general type with $\chi\leq 0$}
\author{Jungkai A. Chen and Meng Chen}

\address{\rm
Department of Mathematics, National Taiwan University, Taipei,
106, Taiwan} \email{jkchen@math.ntu.edu.tw}

\address{\rm Institute of Mathematics, School of Mathematical
Sciences, Fudan University, Shanghai, 200433, People's Republic of
China} \email{mchen@fudan.edu.cn}


\thanks{The first author was partially supported by Taida Institute for Mathematical Sciences,
National Center for Theoretical Sciences, Taipei Office, and
National Science Council of Taiwan. The second author was supported
by Program for New Century Excellent Talents in University
(\#NCET-05-0358), National Outstanding Young Scientist Foundation
(\#10625103) and NNSFC (\#10731030)}
\newcommand{\bC}{{\mathbb C}}
\newcommand{\bQ}{{\mathbb Q}}
\newcommand{\bP}{{\mathbb P}}
\newcommand{\roundup}[1]{\lceil{#1}\rceil}

\newcommand\Alb{\text{\rm Alb}}
\newcommand\alb{\text{\rm alb}}
\newcommand\Pic{\text{\rm Pic}}
\newcommand\Vol{\text{\rm Vol}}

\newcommand\ot{{\otimes}}
\newcommand\OO{{\mathcal{O}}}

\newtheorem{thm}{Theorem}[section]
\newtheorem{lem}[thm]{Lemma}
\newtheorem{cor}[thm]{Corollary}

\theoremstyle{definition}
\newtheorem{defn}[thm]{Definition}
\newtheorem{setup}[thm]{}
\newtheorem{question}[thm]{Question}
\newtheorem{exmp}[thm]{Example}

\newtheorem{rem}[thm]{Remark}
\theoremstyle{remark}

\begin{document}
\begin{abstract} We prove that the canonical volume $K^3\geq
\frac{1}{30}$ for all 3-folds of general type with
$\chi(\OO)\leq 0$. This bound is sharp.
\end{abstract}
\maketitle
\pagestyle{myheadings} \markboth{\hfill J. Chen and M. Chen
\hfill}{\hfill The canonical volume of 3-folds of general type\hfill}

\section{\bf Introduction}
Let $V$ be a nonsingular projective 3-fold of general type.
According to Mori's Minimal Model Program (see for instance
 \cite{KMM, K-M, Reid83}), $V$ has at
least one minimal model $X$ which is normal projective with at
worst ${\bQ}$-factorial terminal singularities. Denote by
$K^3:=K_X^3$. Since it is uniquely determined by the birational
equivalence class of $V$, $K^3$ is usually referred to as {\it the
canonical volume of $V$}, also written as $\Vol(V)$. In the study
of 3-folds of general type, a major difficulty arises when $K^3$
is only a small rational number, rather than an integer. For
example, among known ones by Fletcher-Reid (cf. \cite{C-R}, p151),
$\text{Vol}(V)$ could be as small as $\frac{1}{420}$. It is a fact
that the birational invariant $\Vol(V)$ strongly affects the
geometry of $V$. So  a natural and interesting question is to find
the sharp lower bound $v_3$ of $K^3$ among all those nonsingular
3-folds $V$ of general type.

There have been some relevant known results already:
\begin{itemize}
\item
 There exists a constant $v_3>0$ such that $\Vol(V) \ge v_3$ for all threefolds of general type.
 This is proved by Hacon and
$\text{M}^{\text{c}}$Kernan \cite{H-M}, Takayama \cite{Tak}
and Tsuji \cite{Tsuji};
\item
 It is proved by the second author \cite{mathann}
that $\Vol(V)\geq \frac{1}{3}$ for all 3-folds of general
type with $p_g(V):=\dim H^3(V,\OO_V)\geq 2$ and the bound
``$\frac{1}{3}$'' is sharp.
\end{itemize}

In this paper we would like to prove the following:
\begin{thm}\label{main} Let $V$ be a nonsingular projective
3-fold of general type with $\chi(\OO_V)\leq 0$. Then
\begin{itemize}
\item
[(i)] $\Vol(V)\geq \frac{1}{30}.$
\item
[(ii)] When $\Vol(V)=\frac{1}{30}$, $V$ has the invariants:
$p_g(V)=1$, $q(V)=0$, $\chi(\OO_V)=0$, $P_2(V)=1$,
$P_3(V)=2$, $P_4(V)=3$ and $P_5(V)=4$. Furthermore any
minimal model of $V$ has exactly 3 virtual baskets of
singularities (in the sense of Reid): $1\times
\frac{1}{2}(1,-1,1)$, $1\times \frac{1}{3}(1,-1,1)$,
$1\times\frac{1}{5}(1,-1,1)$.
\end{itemize}
\end{thm}

The next example shows that the lower bound of $\Vol(V)$ in
Theorem \ref{main}(i) is optimal.

\begin{exmp}\label{1/30} (cf. \cite{C-R}, p151) The
canonical hypersurface $X_{28}\subset \bP(1,3,4,5,14)$ has
the canonical volume $K^3=\frac{1}{30}$, $p_g=1$, $q=0$,
$\chi(\OO_{X_{28}})=0$. $X_{28}$ has 3 terminal
singularities: $1\times \frac{1}{2}(1,-1,1)$,
$1\times\frac{1}{3}(1,-1,1)$, $1\times\frac{1}{5}(1,-1,1)$.
\end{exmp}

This note also contains some effective results. For example
we will prove the following:

\begin{cor}\label{effective} Let $V$ be a nonsingular projective
3-fold  of general type with $q:=h^1(\OO_V)>0$. Then
$\Vol(V)\geq \frac{1}{22}$.
\end{cor}

The paper is organized as the following.  In section 2, we study the
pluricanonical maps. We obtain,  in Theorem \ref{p3}, a lower bound
$> \frac{1}{30}$ when the plurigenera are large. In section 3, we
consider irregular threefolds. Combining results obtained in these
two sections, the only unknown case has the information:
$\chi(\OO_X)=0$, $q(X)=0$, $p_g(X)=1$ and $P_5(X) <5$. Thus in
section 4 we classify all possible types of singularities and hence
are able to complete the proof of the main theorem.

Throughout the paper $\sim$ means linear equivalence while
$\equiv$ denotes the numerical one.

\section{\bf Bounding $K^3$ via $\varphi_m$}

In order to get an effective lower bound of $K^3$ we need
to study the $m$-canonical map $\varphi_{m}$. Let $X$ be a
minimal projective 3-fold of general type (admitting at
worst $\bQ$-factorial terminal singularities) with
$$P_{m_0}=P_{m_0}(X):=\dim_{\bC} H^0(X,\OO_X(m_0K_X))\geq 2$$ for some
integer $m_0>0$, where $K_X$ is a canonical divisor of $X$.

\begin{setup}\label{setup}{\bf Set up for $\varphi_{m_0}$.} We
study the $m_0$-canonical map
$\varphi_{m_0}:X\dashrightarrow \bP^{P_{m_0}-1}$ which is
only a rational map. First of all we fix an effective Weil
divisor $K_{m_0}\sim m_0K_X$. By Hironaka's  theorem, we
can take successive blow-ups $\pi: X'\rightarrow X$ such
that:
\begin{itemize}
\item
[(i)] $X'$ is smooth;
\item
[(ii)] the movable part of $|m_0K_{X'}|$ is base point
free;
\item
[(iii)] the support of $\pi^*(K_{m_0})$ is of simple normal
crossings.
\end{itemize}

Set $g_{m_0}:=\varphi_{m_0}\circ\pi$. Then $g_{m_0}$ is a morphism
by assumption. Let $X'\overset{f}\longrightarrow
B\overset{s}\longrightarrow W'$ be the Stein factorization of
$g_{m_0}$ with $W'$ the image of $X'$ through $g_{m_0}$. In
summary, we have the following commutative diagram:\medskip

\begin{picture}(50,80) \put(100,0){$X$} \put(100,60){$X'$}
\put(170,0){$W'$} \put(170,60){$B$}
\put(112,65){\vector(1,0){53}}
\put(106,55){\vector(0,-1){41}}
\put(175,55){\vector(0,-1){43}}
\put(114,58){\vector(1,-1){49}}
\multiput(112,2.6)(5,0){11}{-} \put(162,5){\vector(1,0){4}}
\put(133,70){$f$} \put(180,30){$s$} \put(95,30){$\pi$}
\put(130,-5){$\varphi_{m_0}$}\put(136,40){$g_{m_0}$}
\end{picture}
\bigskip

We recall the definition of $\pi^*(K_X)$ and denote by
$r(X)$ the Cartier index of $X$. Then
$r(X)K_{X'}=\pi^*(r(X)K_X)+E_{\pi}$ where $E_{\pi}$ is a
sum of exceptional divisors. One defines
$\pi^*(K_X):=K_{X'}-\frac{1}{r(X)}E_{\pi}$. So, whenever we
take the round up of $m\pi^*(K_X)$, we always have
$\roundup{m\pi^*(K_X)}\leq mK_{X'}$ for any integer $m>0$.
We may write $m_0K_{X'}=\pi^*(m_0K_X)+E_{
m_0}=M_{m_0}+Z_{m_0},$ where $M_{m_0}$ is the movable part
of $|m_0K_{X'}|$, $Z_{m_0}$ the fixed part and $E_{m_0}$ an
effective ${\mathbb Q}$-divisor which is a ${\mathbb
Q}$-sum of distinct exceptional divisors. We may also write
$m_0\pi^*(K_X)= M_{m_0}+E_{m_0}',$ where
$E_{m_0}'=Z_{m_0}-E_{m_0}$ is an effective ${\mathbb
Q}$-divisor.

If $\dim(B)\geq 2$, a general member $S$ of $|M_{m_0}|$ is a
nonsingular projective surface of general type by Bertini's
theorem and by the easy addition formula for Kodaira dimension.

If $\dim(B)=1$, a general fiber $S$ of $f$ is an
irreducible smooth projective surface of general type,
still by the easy addition formula for Kodaira dimension.
We may write
$$M_{m_0}=\underset{i=1}{\overset{a_{m_0}}\sum}S_i\equiv
a_{m_0}S$$ where the $S_i$ is a smooth fiber of $f$ for all
$i$ and $a_{m_0}\ge P_{m_0}(X)-1$.

In both cases we call $S$ {\it a generic irreducible
element of} $|M_{m_0}|$. Denote by $\sigma:
S\longrightarrow S_0$ the blow-down onto the smooth minimal
model $S_0$.
\end{setup}

\begin{setup}\label{p,b}{\bf Assumptions}. We need some
assumptions to estimate $K^3$.
\begin{itemize}
\item
[(1)] Keep the same notations as above, we define
$$p=\begin{cases}   1 &\text{if}\ \dim (B)\geq 2\\
 a_{m_0} &\text{if}\ \dim(B)=1.
 \end{cases}$$

\item
[(2)] Take a generic irreducible element $S$ of
$|M_{m_0}|$. Assume that $|G|$ is a movable complete linear
system on $S$. Also assume that a generic irreducible
element $C$ of $|G|$ is smooth.
\item
[(3)] Assume there is a rational number $\beta>0$ such that
$\pi^*(K_X)|_S-\beta C$ is numerically equivalent to an
effective $\bQ$-divisor on $S$.
\end{itemize}

Set $\alpha=(m-1-\frac{m_0}{p}-\frac{1}{\beta})\xi$ and
$\alpha_0:=\roundup{\alpha}$.
\medskip

Under Assumptions \ref{p,b}, one has
$$K^3\geq \frac{p}{m_0}\pi^*(K_X)^2\cdot S\geq
\frac{p\beta}{m_0}(\pi^*(K_X)\cdot C). \eqno{(2.1)}$$

So it suffices to estimate the rational number
$\xi:=(\pi^*(K_X)\cdot C)_{X'}$.
\end{setup}

We need the following theorem to study the lower bound of
$\xi$:

\begin{thm}\label{technical} Let $m>0$ be an integer. Under
Assumptions \ref{p,b}, the inequality
$$\xi\geq \frac{\deg(K_C)+\alpha_0}{m}$$
holds if one of the following conditions is satisfied:
\begin{itemize}
\item
[(i)] $\alpha_0\geq 2$;
\item
[(ii)] $\alpha>0$ and $C$ is an even divisor on $S$.
\end{itemize}
\end{thm}
\begin{proof} We consider the sub-linear system
$$|K_{X'}+\roundup{(m-1)\pi^*(K_X)-\frac{1}{p}E_{m_0}'}|\subset
|mK_{X'}|.$$ Take a generic irreducible element $S$ of
$|M_{m_0}|$. Noting that
$(m-1)\pi^*(K_X)-\frac{1}{p}E_{m_0}'-S\equiv
(m-1-\frac{m_0}{p})\pi^*(K_X)$ is nef and big whenever $\alpha>0$,
the Kawamata-Viehweg vanishing theorem \cite{V,KV} gives the
surjective map
$$H^0(X',K_{X'}+\roundup{(m-1)\pi^*(K_X)-\frac{1}{p}E_{m_0}'})$$
$$\longrightarrow  H^0(S,
K_{S}+\roundup{(m-1)\pi^*(K_X)-S-\frac{1}{p}E_{m_0}'}|_{S}).\eqno
(2.2)$$

Now consider a generic irreducible element $C\in |G|$. By
assumption there is an effective ${\mathbb Q}$-divisor $H$
on $S$ such that
$$\frac{1}{\beta}\pi^*(K_X)|_{S}\equiv C+H.$$
By the vanishing theorem again, whenever
$m-1-\frac{m_0}{p}-\frac{1}{\beta}>0$ which yields that
$$((m-1)\pi^*(K_X)-S-\frac{1}{p}E_{m_0}')|_{S}-C-H\equiv (m-1-
\frac{m_0}{p}-\frac{1}{\beta})\pi^*(K_X)|_{S}$$ is nef and
big, we have the surjective map
$$H^0(S,
K_{S}+\roundup{((m-1)\pi^*(K_X)-S-\frac{1}{p}E_{m_0}')|_{S}-H})$$
$$\longrightarrow  H^0(C, K_C + D)\hskip6.8cm \eqno (2.3)$$
where
$D:=\roundup{((m-1)\pi^*(K_X)-S-\frac{1}{p}E_{m_0}')|_{S}-C-H}|_C$
is a divisor on $C$. Noting that $C$ is nef on $S$, we have
$\deg(D)\geq\alpha$ and thus $\deg(D)\geq \alpha_0$.

Whenever either $\deg(D)\geq 2$ or $C$ is an even divisor
and $m-1-\frac{m_0}{p}-\frac{1}{\beta}>0$ ($\deg(D)\geq 2$
automatically follows), $|K_C+D|$ is base point free by the
curve theory. Denote by $|N_m|$ the movable part of
$|K_{S}+\roundup{((m-1)\pi^*(K_X)-S-\frac{1}{p}E_{m_0}')|_{S}-H}|$.
Applying Lemma 2.7 of \cite{MPCPS} to surjective maps (2.2)
and (2.3), one has
$$m\pi^*(K_X)|_{S}\ge N_m\ \ \text{and}\ \ (N_m\cdot C)_{S}\ge 2g(C)-2
+\deg (D).$$ So $m\xi\geq \deg(K_C)+\alpha_0$. We are done.
\end{proof}

\begin{rem}\label{beta} A technical problem in utilizing
Theorem \ref{technical} is to verify Assumptions \ref{p,b}.
To avoid unnecessary redundancy, we only copy several
technical results here without proof. Note that the most
complicated situation is the one with $\dim (B)=1$, in
which case we set $b:=g(B)$, the geometric genus.
\begin{itemize}
\item
[(1)] Usually we will take $G\leq \sigma^*(K_{S_0})$
whenever $p_g(S)\geq 2$ or $G=2\sigma^*(K_{S_0})$,
$4\sigma^*(K_{S_0})$ otherwise;
\item
[(2)] When $g(B)>0$, it is proved in Lemma 3.4 of
\cite{Chen-Zuo} that $$\pi^*(K_X)|_S\sim
\sigma^*(K_{S_0}).$$ So one may take $\beta=1$ or
$\frac{1}{2}$ or $\frac{1}{4}$.
\item
[(3)] When $g(B)=0$, it is proved in Lemma 3.4 of
\cite{Chen-Zuo} that
$$\pi^*(K_X)|_S-\widetilde{\beta}_n\sigma^*(K_{S_0})$$ is
numerically equivalent to an effective $\bQ$-divisor for a
sequence of positive rational numbers
$\{\widetilde{\beta}_n\}$ with $\widetilde{\beta}_n\mapsto
\frac{p}{m_0+p}$. So one may take
$\beta=\widetilde{\beta}_n$ or
$\frac{1}{2}\widetilde{\beta}_n$ or
$\frac{1}{4}\widetilde{\beta}_n$ accordingly.
\end{itemize}
\medskip

(\ddag) Note that a special situation (with $m_0=1$) of
this theory has already appeared in \cite{mathann}, as 2.8
and Lemma 3.4.
\end{rem}

Now we are ready to estimate $\xi$ and $K^3$.

\begin{thm}\label{p3} Let $V$ be a nonsingular projective
3-fold of general type. Then
\begin{itemize}
\item
[(i)] $\Vol(V)\geq \frac{1}{22}$ if $P_4(V)\geq 5$;
\item
[(ii)] $\Vol(V)\geq \frac{1}{25}$ if $P_5(V)\geq 5$,
$p_g(V)>0$ and $\dim(B)\geq 2$.
\item
[(iii)] $\Vol(V)\geq \frac{8}{45}$ if $P_5(V)\geq 5$,
$p_g(V)>0$ and $\dim(B)=1$.
\end{itemize}
\end{thm}
\begin{proof} Take a minimal model $X$ of $V$. We study $|mK|$ on $X$.
Keep the same set up as in \ref{setup}.
\medskip

\noindent{\bf Part (i)}. We can take $m_0=4$. We will study
according to the value of $\dim (B)$. Take a generic
irreducible element $S$ of $|M_4|$.

If $\dim (B)=3$, we know that $p=1$ by definition. In this
case we know $S\sim M_5$ and that $|S|$ gives a generically
finite morphism. Set $G:=S|_S$. Then $|G|$ is base point
free and $\varphi_{|G|}$ gives a generically finite map. So
a generic irreducible element $C$ of $|G|$ is a smooth
curve of genus $\geq 2$. If $\varphi_{|G|}$ gives a
birational map, then $\dim \varphi_{|G|}(C)=1$ for a
general member $C$. The Riemann-Roch and Clifford's theorem
on $C$ says $C^2=G\cdot C\geq 2$. If $\varphi_{|G|}$ gives
a generically finite map of degree $\geq 2$, since
$h^0(S,G)\geq h^0(X',S)-1\geq 4$, Lemma 2.2 of \cite{JMSJ3}
gives $C^2\geq 2h^0(S,G)-4\geq 4$. Anyway we have $C^2\geq
2$. So $\deg(K_C)=(K_S+C)C>2C^2\geq 4$. We see
$\deg(K_C)\geq 6$ because it is even. One may take
$\beta=\frac{1}{4}$ since $4\pi^*(K_X)|_S\geq C$. Now if we
take a very big $m$ such that $\alpha>1$ then Theorem
\ref{technical} gives:
$$m\xi\geq \deg(K_C)+(m-1-m_0-\frac{1}{\beta})\xi.$$
This gives $\xi\geq \frac{2}{3}$. If we take $m=11$. Then
$\alpha=2\xi>1$. Theorem \ref{technical} says $\xi\geq
\frac{8}{11}$. So inequality (2.1) gives $K^3\geq \frac{1}{22}$.

If $\dim(B)=2$, we know that $|G|:=|S|_S|$ is composed with
a pencil of curves. A generic irreducible element $C$ of
$|G|$ is a smooth curve of genus $\geq 2$, so
$\deg(K_C)\geq 2$. Furthermore we have $h^0(S,G)\geq
h^0(X',S)-1\geq 4$. So $G\equiv \widetilde{a} C$ for
$\widetilde{a}\geq h^0(S,G)-1\geq 3$. This means
$4\pi^*(K_X)|_S\geq S|_S\geq_{\text{numerically}} 3C$. So
we may take $\beta=\frac{3}{4}$. Now take a very big $m$.
Theorem \ref{technical} gives $\xi\geq \frac{6}{19}$. Take
$m=10$. Then $\alpha\geq \frac{22}{19}>1$. We get $\xi\geq
\frac{2}{5}$. So inequality (2.1) gives $K^3\geq
\frac{3}{40}>\frac{1}{22}$.
\medskip

If $\dim (B)=1$ and $b=g(B)>0$, there is an induced
fibration $f:X'\longrightarrow B$. Recall that $S$ is a
general fiber of $f$ and $S$ can be a nonsingular surface
of general type of any numerical type. One has $p=a_4\geq
P_4\geq 5$ by the Riemann-Roch and Clifford's theorem. Set
$G:=4\sigma^*(K_{S_0})$. Because $|4K_{S_0}|$ is base point
free by Bombieri \cite{Bom}, $|G|=|4\sigma^*(K_{S_0})|$ is
also base point free. Denote by $C$ a generic irreducible
element of $|G|$. Then $C$ is smooth and
$\deg(K_C)=(K_S+C)C\geq
(\pi^*(K_X)|_S+C)C=(\sigma^*(K_{S_0})+C)C=20\sigma^*(K_{S_0})^2\geq
20$. By Remark \ref{beta}(2), we can take
$\beta=\frac{1}{4}$ since $\pi^*(K_X)|_S\sim
\frac{1}{4}(4\sigma^*(K_{S_0}))$. Now if we take a very big
$m$, Theorem \ref{technical} gives $\xi\geq
\frac{100}{29}$. Inequality (2.1) gives $K^3\geq
\frac{125}{116}>\frac{1}{22}$.
\medskip

If $\dim (B)=1$ and $b=g(B)=0$, we have $p\geq P_5-1\geq
4$. Take $m_0=4$.
%
%
We set
$G:=4\sigma^*(K_{S_0})$. The surface theory (see e.g. Bombieri
\cite{Bom}, or Reider \cite{Reider}) tells us that
$|G|$ is base point free and a generic irreducible element
$C$ of $|G|$ is a smooth curve. Because
$$\deg(K_C)=(K_S+C)C\geq (\pi^*(K_X)|_S+C)C>C^2\geq 16,$$
 we see $\deg(K_C)\geq 18$ due to the evenness.
We know that
$\pi^*(K_X)|_S-\widetilde{\beta}_n\sigma^*(K_{S_0})$ is
numerically equivalent to an effective $\bQ$-divisor for a
rational number sequence $\{\widetilde{\beta}_n\}$ with
$\widetilde{\beta}_n\mapsto \frac{p}{p+4}\geq \frac{1}{2}$.
Take $\beta_n:=\frac{1}{4}\widetilde{\beta}_n$. Then
$\pi^*(K_X)|_S-\beta_n C$ is numerically equivalent to an
effective $\bQ$-divisor. We know $\beta_n\mapsto
\frac{1}{8}$ whenever $p=4$. When $m$ is very big, Theorem
\ref{technical} gives $\xi\geq \frac{9}{5}$. Take $m=11$.
Then $\alpha\geq \xi>1$. One gets $\xi\geq
\frac{20}{11}>\frac{9}{5}$. So inequality (2.1) gives
$K^3\geq \frac{5}{22}>\frac{1}{22}$.

Comparing what we have proved, we see $K^3\geq
\frac{1}{22}$.
\medskip

\noindent {\bf Part (ii)}. Take $m_0=5$. We have $p=1$ by
definition. A general member $S\in |M_5|$ is a nonsingular
projective surface of general type. Set $G:=S|_S$.

If $\varphi_5$ is generically finite, then $\varphi_{|G|}$
is either birational or generically finite of degree $\geq
2$. We have the following argument:

\begin{quote}
($\sharp$)  If $\varphi_{|G|}$ gives a birational map, then
clearly $h^0(S,G)\geq 4$ because $S$ is of general type.
Since $p_g(V)>0$, we know $G\leq (K_{X'}+S)|_S=K_S$. And
because $G$ is nef, Lemma 2.1 of \cite{JMSJ3} says $C^2\geq
3h^0(S,G)-7\geq 5$. When $|S|$ gives a generically finite
map of degree $\geq 2$, then Lemma 2.2 of \cite{JMSJ3}
gives $C^2\geq 2h^0(S,G)-4$.
\end{quote}

Because $h^0(S,G)\geq h^0(X',S)-1\geq 4$, the argument
($\sharp$) says $C^2\geq 4$. We get
$\deg(K_C)=(K_S+C)C>2C^2\geq 8$ noting that $K_{X'}|_S\cdot
C\geq \pi^*(K_X)|_S\cdot C>0$ by the Hodge Index Theorem.
Actually we have $\deg(K_C)\geq 10$ since it is even. On
the other hand, we can take $\beta=\frac{1}{5}$ since
$5\pi^*(K_X)|_S\geq C$. Now take a very big $m$. Theorem
\ref{technical} gives $\xi\geq \frac{10}{11}$. Take $m=13$.
Then $\alpha\geq \frac{20}{11}>1$. We get $\xi\geq
\frac{12}{13}$. So inequality (2.1) gives $\xi\geq
\frac{12}{25\cdot 13}$. In fact a similar calculation says
$\xi\geq \frac{l}{l+1}$ for all $l\geq 12$. Thus $\xi\geq
1$ and (2.1) gives $K^3\geq \frac{1}{25}$.

If $\dim(B)=2$, we know that $|G|:=|S|_S|$ is composed with
a pencil of curves. A generic irreducible element $C$ of
$|G|$ is a smooth curve of genus $\geq 2$, so
$\deg(K_C)\geq 2$. Furthermore we have $h^0(S,G)\geq
h^0(X',S)-1\geq 4$. So $G\equiv \widetilde{a} C$ for
$\widetilde{a}\geq h^0(S,G)-1\geq 3$. This means
$5\pi^*(K_X)|_S\geq S|_S\geq_{\text{numerically}} 3C$. So
we may take $\beta=\frac{3}{5}$. Now take a very big $m$.
Theorem \ref{technical} gives $\xi\geq \frac{6}{23}$. Take
$m=12$. Then $\alpha\geq \frac{26}{23}>1$. We get $\xi\geq
\frac{1}{3}$. Take $m=11$. Then $\alpha\geq
\frac{10}{9}>1$. We get $\xi\geq \frac{4}{11}$. So
inequality (2.1) gives $K^3\geq \frac{1}{25}\cdot
\frac{12}{11}>\frac{1}{25}$.
\medskip

\noindent {\bf Part (iii)}. Take $m_0=5$. Parallel to the
last parts in the proof of (i), we can discuss according to
the value of $b=g(B)$. So we have more or less a redundant
calculation as follows.

If $b=g(B)>0$, there is an induced fibration $f:X'\longrightarrow
B$. Because $p_g(V)>0$, one sees $p_g(S)>0$. One has $p=a_5\geq
P_5\geq 5$ by the Riemann-Roch and Clifford's theorem. Set
$G:=2\sigma^*(K_{S_0})$. Because $|2K_{S_0}|$ is base point free
(see Theorem 3.1 in \cite{Ci}), $|G|=|2\sigma^*(K_{S_0})|$ is also
base point free. Denote by $C$ a generic irreducible element of
$|G|$. Then $C$ is smooth and $\deg(K_C)=(K_S+C)C>
4\sigma^*(K_{S_0})^2\geq 4$. So actually $\deg(K_C)\geq 6$. By
Remark \ref{beta}(2), we can take $\beta=\frac{1}{2}$ since
$\pi^*(K_X)|_S\sim \frac{1}{2}(2\sigma^*(K_{S_0}))$. Now if we take
a very big $m$, Theorem \ref{technical} gives $\xi\geq \frac{3}{2}$.
Inequality (2.1) gives $K^3\geq \frac{3}{4}$.

If $b=g(B)=0$, we have $p\geq P_5-1\geq 4$. Take $m_0=5$.
Because $p_g(V)>0$, one sees $p_g(S)>0$. We still set
$G:=2\sigma^*(K_{S_0})$. We have $\deg(K_C)\geq 6$.
Similarly we only have to find a suitable $\beta$. Remark
\ref{beta}(3) says that one can find a sequence of positive
rational numbers $\{\widetilde{\beta}_n\}$ with
$\widetilde{\beta}_n\mapsto \frac{p}{p+5}\geq \frac{4}{9}$
such that
$\pi^*(K_X)|_S-\widetilde{\beta}_n\sigma^*(K_{S_0})$ is
numerically equivalent to an effective $\bQ$-divisor. Set
$\beta_n:=\frac{1}{2}\widetilde{\beta}_n$. Then
$\pi^*(K_X)|_S-\beta_n C$ is numerically equivalent to an
effective $\bQ$-divisor. We know $\beta_n\mapsto
\frac{2}{9}$ whenever $p=4$. When $m$ is very big, Theorem
\ref{technical} gives $\xi\geq \frac{8}{9}$. Take $m=8$.
Then $\alpha\geq \frac{10}{9}>1$. We get $\xi\geq 1$. So
inequality (2.1) gives $K^3\geq \frac{8}{45}$. This
completes the proof.
\end{proof}

\begin{rem} One may remove extra condition: $p_g(V)>0$
in Theorem \ref{p3} (ii) and (iii) to obtain parallel, but weaker
results. We omit the details simply because it is not used in the
proof of the main theorem.
\end{rem}

\section{\bf  Irregular 3-folds of general type}

In this section, we study the canonical volume of irregular
threefolds of general type. Let $X$ be a nonsingular
projective threefold of general type and $a:X \to \alb(X)$
the Albanese map. By running the minimal model program, one
easily see that the Albanese map factors through its
minimal model. So we may and do assume that $X$ is a
minimal ($K_X$ nef) threefold of general type with
$\bQ$-factorial terminal singularities.


In the study of pluricanonical systems on irregular
threefolds, the most unpleasant case is when the Albanese
map $a: X \to \Alb(X)$ is surjective onto an elliptic curve
$E$ with general fiber $F$ of type $(K_F^2, p_g(F))=(1,2)$.

\begin{thm}\label{1/9} Let $X$ be a minimal
3-fold of general type with $q(X)=1$ and the general fiber
of $a:X \to \Alb(X)$ is of $(1,2)$ type. Then the canonical
volume $\text{Vol}(X)=K_X^3\geq \frac{1}{9}$.
\end{thm}

Before proving the main result, we would like to recall
some notion and results in \cite{JC-H-mz}.

\begin{defn}
For any vector bundle $E$ on an elliptic curve, we write $E
= \oplus E_i$ for its decomposition into indecomposable
vector bundles. We define $\nu(E):=\min \mu(E_i)$, where
$\mu(E_i)=\frac{deg(E_i)}{rk(E_i)}$.
\end{defn}

\begin{lem}\label{mucomp}(\cite{JC-H-mz}, Lemma 4.8)
Let $E_1, E_2$ be indecomposable vector bundles on an
elliptic curve. If ${\rm Hom}(E_1,E_2) \neq 0$, then
$\mu(E_2) \ge \mu (E_1)$. In particular $\nu(E_2) \ge \nu
(E_1)$ if $E_1 \to E_2$ is a surjective map of vector
bundles.
\end{lem}

\begin{defn} A coherent sheaf $\mathcal {F}$ on an abelian
variety $A$ is said to be $IT^0$ if $H^i(A,\mathcal
{F}\otimes P)=0$ for all $i>0$ and all $P\in
\text{Pic}^0(A)$.
\end{defn}

\begin{lem}\label{numin} (\cite{JC-H-mz}, Lemma 4.10)
Let $E$ be an $IT^0$ vector bundle on an elliptic curve
which admits a short exact sequence $$ 0 \to F \to E \to Q
\to 0$$ of coherent sheaves such that $Q$ has generic rank
$\le 1$. Then
 $\nu(E) \ge \min\{1, \nu(F)\}$.
\end{lem}

\begin{setup}{\bf Multiplication maps $\varphi_{m,n}$ and
$\psi_{m,n}$.} Let $R_m:= H^0(F, \omega_F^m)$ and $E_m:=a_*
\omega_X^m$. By Lemma 4.1 of \cite{JC-H-mz}, $E_m$ is an
$IT^0$ vector bundle of rank $P_m(F)$ for all $m\geq 2$. We
also remark that $\nu(E_1) \ge 0$ by the semipositivity
theorem (see Viehweg \cite{V1}) and Atiyah's description of
vector bundles over elliptic curves (cf. \cite{At}). We
consider the multiplication map of pluricanonical systems
on fibers
$$\varphi_{m,n}: R_m \ot R_n \to R_{m+n}.$$
This induces a map
$$\psi_{m,n}: E_m \ot E_n \to E_{m+n}.$$
Clearly if cokernel of $\varphi_{m,n}$ has dimension $\leq
r$, then cokernel of $\psi_{m,n}$ has rank $\leq r$.
\end{setup}

\begin{setup}{\bf Surfaces of (1,2) type.} Let $F$ be a nonsingular
minimal projective surface of general type with
$(K_F^2,p_g(F))=(1,2)$. It's well-known that $|K_F|$ has
only one base point $z$ and $|2K_F|$ is base point free
(cf. \cite{Bom}).
\end{setup}

We recall the following result in \S 5 of \cite{JC-H-mz}.

\begin{lem}\label{phi1}(\cite{JC-H-mz}, p353, line -6)
Assume that a general fiber $F$ of the fibration $a: X \to
\Alb(X)$ is a surface of (1,2) type. Then
 $ \varphi_{1,m-1}: R_1 \ot R_ {m-1}\mapsto R_1R_ {m-1} \subset R_{m}$ has
codimension $ \le 1$ for all $m \geq 1 $ and
$\varphi_{1,2}$ is surjective.
\end{lem}

($\dag$) Clearly $R_1 R_2 =R_3$ implies that $R_2 R_2$,
which contains $R_1 R_1 R_2=R_1R_3$, has codimension $\leq
1$ in $R_4$.

Moreover, we have the following:

\begin{lem}\label{phi2} Assume that a general fiber $F$ of
the fibration $a: X \to \Alb(X)$ is a surface of (1,2) type. Then
the multiplication map $\varphi_{2,m-2}: R_2 \ot R_ {m-2} \mapsto
R_{m}$ is surjective whenever $m \ge 8$ and has codimension $\leq
1$ if $m =7$.
\end{lem}

\begin{proof}
We follow Bombieri's argument on projective normality.

Fix two sections $s, s_1 \in R_2$ with smooth curves
$C:=\text{div}(s)$, $C_1:=\text{div}(s_1)$ and assume that
$Z:=div(s_1) \cap C$ consists of $4$ distinct points, we
have the exact sequences:
$$ 0 \to R_{m-2} \stackrel{s}{\to} R_m \stackrel{r_C}{\to} H^0(C,
mK_{F}|_C) \to 0,\eqno{(3.1)}$$
$$ 0 \to R_{m-2} \stackrel{s_1}{\to} R_m \stackrel{r_{C_1}}{\to} H^0(C_1,
mK_{F}|_{C_1}) \to 0,\eqno{(3.2)}$$
$$ 0 \to R_{m-4} \stackrel{s}{\to} R_{m-2} \stackrel{r_C}{\to} H^0(C,
(m-2)K_{F}|_C) \to 0, \eqno{(3.3)}$$ thanks to the
vanishing of $H^1(F, tK_F)$ for all $t \ge 0$.
We also have the following exact sequence:
$$ 0 \to H^0(C,(m-2)K_{F}|_C) \stackrel{\widetilde{s}_1}{\to} H^0(C,
mK_{F}|_C) \stackrel{r_Z}{\to} H^0(Z,\OO_Z).\eqno{(3.4)}$$
Among the four exact sequences, one can find the
commutative relation: $\widetilde{s}_1\circ r_C=r_C\circ
s_1$.

$$\begin{CD}
R_{m-2}@> s_1>>R_m\\
@V r_C VV @VV r_C V\\
H^0(C,(m-2)K_F|_C)@>> \widetilde{s}_1 >H^0(C,mK_F|_C)\\
\end{CD}$$

One knows $\widetilde{s}_1\circ
r_C(R_{m-2})=\widetilde{s}_1 H^0(C,(m-2)K_F|_C)$ has
codimension $\leq 4$ in the space $H^0(C, mK_F|_C)$ since
$\dim H^0(Z,\OO_Z)=4$. So $r_C\circ s_1(R_{m-2})$ has
codimension $d'\leq 4$ in $H^0(C, mK_F|_C)$. Since
\begin{eqnarray*}
d'&=&h^0(C,mK_F|_C)-\dim r_C(s_1R_{m-2})\\
&\geq& \dim R_m-\dim s R_{m-2}-\dim s_1 R_{m-2}.
\end{eqnarray*}
It follows that $sR_{m-2} + s_1 R_{m-2} \subset R_m$ has
codimension $\leq 4$.

Moreover, we consider
\begin{eqnarray*} 0 &\to&
H^0(C,(m-4)K_{F}|_C) \stackrel{\widetilde{s}_1}{\to} H^0(C,
(m-2)K_{F}|_C)\\
 & \stackrel{r_Z}{\to}& H^0(Z,\OO_Z) \to
H^1(C, (m-4)K_{F}|_C).
\end{eqnarray*}
Since $3K_{F|C}=K_C$ and $\deg(K_{F}|_C)=2$, one sees that
$H^1(C, (m-4)K_{F}|_C)=0$ if $m\geq 8$. When $m=7$,
$H^1(C,3K_{F}|_C)=H^1(C, K_C)$ is one dimensional.
We can take a section $s_2 \in R_2$ such that $s_2$ never
vanishing on $Z$. Set $J=div(s_2) \cap C$ which can be a
union of 4 distinct points. As we have seen the map $r_J:
r_C(s_2 R_{m-2}) \to H^0(J,\OO_J)$ is either surjective
when $m \ge 8$ or having codimension $\le 1$ when $m=7$.
Together with surjectivity of $r_C$, we see that $r_C(s_2
R_{m-2})=\widetilde{s}_2(H^0(C,(m-2)K_F|_C))$ where
$$\widetilde{s}_2: H^0(C, (m-2)K_{F}|_C)\mapsto H^0(C,
mK_F|_C)$$ is defined by the multiplication of $s_2$. This
already means that
$$s R_{m-2}+s_1 R_{m-2}+s_2 R_{m-2}
\subset R_m$$ has codimension $0$ or $\le 1$ if $m \ge 8$
or $=7$ respectively. We are done.
\end{proof}

Now we prove Theorem \ref{1/9}.
\begin{proof}[\bf Proof of Theorem \ref{1/9}]
First of all, $E_2$ is an $IT^0$ vector bundle of rank $4$.
So one has $\nu(E_2) \geq \frac{1}{4}$.

Consider the induced multiplication map $\psi_{2,2}:E_2 \ot
E_2 \to E_4$. Since $\varphi_{2,2}$ has image of
codimension $\leq 1$ by $(\dag)$, it follows that
$\psi_{2,2}$ has cokernel of rank $\leq 1$. We consider the
exact sequence
$$ 0 \to \text{Im}(\psi_{2,2}) \to E_4 \to \text{Coker}(\psi_{2,2}) \to 0.$$

By Lemma \ref{numin} and Lemma \ref{mucomp}, one has
\begin{eqnarray*}
\mu(E_4)& \ge& \min\{ \nu(\text{Im}( \psi_{2,2})),1\} \ge
\min\{\nu(E_2 \ot E_2),1\}\\
& =& \min\{2\nu(E_2),1\} \ge \frac{1}{2}.
\end{eqnarray*}

Next, similarly, we consider $\psi_{4,1}$, then we see that
$$\nu(E_5) \ge  \min\{\nu(E_4) + \nu(E_1), 1\}  \geq \frac{1}{2}.$$
We now consider $\psi_{5,2}$ and the induced exact sequence
$$ 0 \to \text{Im}(\psi_{5,2}) \to E_7 \to \text{Coker}(\psi_{5,2}) \to 0.$$
By Lemma 3.9, one sees that $\text{Coker}(\psi_{5,2})$ has generic rank $\le 1$. Therefore
Lemma 3.5 can be applied. Argue as in the consideration of $\psi_{2,2}$,
we see that $\nu(E_7) \geq
\min\{\nu(E_5) + \nu(E_2),1\} \ge \frac{3}{4}$.

Finally, we consider $\psi_{7,2}$, then similarly we have $\nu(E_9)
\geq 1$.

Now $\nu(E_9) \ge 1$ implies that there is a line bundle
$L$ of degree $1$ with an injection $L \to E_9$. In
particular, $H^0(X, 9K_X \ot f^*L^\vee)$ has a section.
Thus $9K_X \geq  F$. Hence
$$9K_X^3
\geq K_X^2 \cdot F = K_F^2  =1.$$ This completes the proof.
\end{proof}




Combining Theorem \ref{1/9}, results in \cite{JC-H-mz}, and
Theorem \ref{p3}, we are able to get a lower bound of the
canonical volume for all those irregular threefolds.

\begin{cor}\label{q>0} Let $V$ be a nonsingular projective
irregular 3-fold  of general type. Then $\Vol(V)\geq
\frac{1}{22}$.
\end{cor}

\begin{proof} We consider 3-folds of general type with $q(V) >0$.
Then there is a non-trivial Albanese map $a: V \to
\Alb(V)$. If the general fiber has dimension $\leq 1$, then
by Proposition 2.9 of \cite{JC-H-mz}, $|4K_V+P|$ is
birational for general $P \in \Pic^0(V)$. In particular,
$h^0(V, \OO_V(4K_V) \ot P) \ge 4$. However, it's in fact
$\ge 5$ because otherwise it gives a birational map onto
$\mathbb{P}^3$, which is not of general type. By the upper
semicontinuity of cohomology, we have $P_4(V)=h^0(V,
\OO(4K_V)) \ge 5$. Now by Theorem \ref{p3} (i), we get
$\Vol(V) \ge \frac{1}{22}$.

We now assume that  the Albanese map has $1$-dimensional
image. Let $f: V \to  H$ be an induced fibration from the
Stein factorization of $a$. We now consider the case that
$g(H) \ge 2$. We can take the relative minimal model of
$f$, say $h:X\to H$. So $X$ is birational to $V$. By
Theorem 1.4 of \cite{ohno}, $K_{X/H}:=K_X-h^*(K_H)$ is nef.
In particular $K_X$ is nef and $X$ is minimal. Because
$g(H)\geq 2$, we see that $K_X-2F$ is nef where $F$ is a
general fiber of $h$. So $\Vol(V)=K_X^3\geq 2K_F^2\geq
2>\frac{1}{22}$.

Finally, we consider the case that $g(H)=1$. We remark that
$g(H)=1$ if and only if $q(V)=1$ because if $q(V) \ge 2$,
then either its Albanese image has dimension $\geq 2$ or is
a curve of genus $\geq 2$.

 If $F$ is not of the type $(1,2)$, then $|4K_F|$ is
birational according to Bombieri's classification. By
Theorem 2.8 of \cite{JC-H-mz}, $|4K_V+P|$ is birational for
general $P \in \Pic^0(X)$. So we get $\Vol(V) \geq
\frac{1}{22}$ as above.

It remains to consider the case that $F$ is of type
$(1.2)$. By Theorem \ref{1/9}, we have $\Vol(V) \ge
\frac{1}{9}>\frac{1}{22}$.
\end{proof}

\section{\bf The case $P_5<5$}

\begin{setup}\label{RR}
First let us recall Reid's plurigenera formula (cf.  \cite{YPG},
p413) for a minimal 3-fold $X$ of general type (with
$\bQ$-factorial terminal singularities):
$$P_m(X)=\frac{1}{12}m(m-1)(2m-1)K_X^3-(2m-1)\chi(\mathcal
{O}_X)+l(m)  \eqno (4.1)$$ where $m$ is an integer $>1$. The
correction term is
$$l(m):=\sum_{Q}l_Q(m):=\sum_{Q}\sum_{j=1}^{m-1}
\frac{\overline{bj}(r-\overline{bj})}{2r},$$ where the sum
$\sum_{Q}$ runs through all baskets Q of singularities of type
$\frac{1}{r}(a,-a,1)$ with the integer $a$ coprime to $r$,
$0<a<r$, $0<b<r$, $ab\equiv 1$ (\text{mod} $r$), $\overline{bj}$
the smallest residue of $bj$ \text{mod} $r$. One can see easily
that $(b,r)=1$. Note by definition that the singularity
$\frac{1}{r}(a,-a,1)$ is a terminal quotient one obtained by a
cyclic group action on $(\bC^3, (0,0,0))$:
$$\varepsilon (x,y,z)=(\varepsilon^a x,\varepsilon^{-a} y, z)$$
where $\varepsilon$ is a fixed $r$-th primitive root of $1$.
Reid's Theorem 10.2 in \cite{YPG} says that the above baskets
$\{\text{Q}\}$ of singularities are in fact virtual (!) and that
one need not worry about the authentic type of all those terminal
singularities on $X$, though $X$ may have non-quotient terminal
singularities. Iano-Fletcher \cite{Fletcher} has shown that the
set of baskets $\{\text{Q}\}$ in Reid's formula is uniquely
determined by $X$.
\end{setup}

In the next context we will always study those 3-folds $X$ with
the following conditions:\medskip

\begin{quote} (*) $p_g=1$, $\chi(\OO)=0$ and $P_5\leq 4$.
\end{quote}

\begin{setup}\label{n} Reid's formula (4.1) tells:
$$P_5>P_4>P_3>P_2>0$$
whenever $\chi(\OO)=0$. So one gets $P_2=1$, $P_3=2$, $P_4=3$ and
$P_5=4$.

We shall classify those $X$ satisfying the condition that
$\chi(\OO_X)=0, P_2(X)=1,P_3(X)=2$.

Recall the plurigenera formula that
$$P_m(X)=\frac{1}{12}m(m-1)(2m-1)K_X^3+ \sum_{Q}\sum_{j=1}^{m-1}
\frac{\overline{b_Qj}(r_Q-\overline{b_Qj})}{2r_Q} .$$

We introduce $$b'_Q :=\left\{ \begin{array}{ll} b_Q, &{\text{if}}
\quad b_Q \le
\frac{1}{2} r_Q; \\
r_Q-b_Q,& {\text{if}}  \quad b_Q > \frac{1}{2} r_Q. \end{array}
\right .$$ Then it's easy to see that
${\overline{b_Qj}(r_Q-\overline{b_Qj})} = {{b'_Qj}(r_Q-{b'_Qj})}$
for $j=1,2$.

For $m=2,3$, we have
$$ 1= P_2(X)= \frac{1}{2}K_X^3 +
\sum_Q\frac{{b'_Q}(r_Q-{b'_Q})}{2r_Q}= \frac{1}{2}K_X^3 +
\frac{1}{2}\sum_Q b'_Q-\frac{1}{2} \sum_Q \frac{{b'}_Q^2}{r_Q},$$
$$ 2= P_3(X) = \frac{5}{2}K_X^3 +
\frac{3}{2}\sum_Q b'_Q-\frac{5}{2} \sum_Q \frac{{b'}_Q^2}{r_Q}.$$

By solving these, we get $$\sum_Q b'_Q=3, \quad \sum_Q
\frac{{b'}_Q^2}{r_Q}=1+K_X^3. \eqno{(4.1)}$$


Moreover, the inequality
$$1=P_2(X)\geq \frac{1}{2}K^3+\sum \frac{r-1}{2r}>
\frac{n}{4}$$ implies $n<4$, where $n$ denotes the number of
baskets. Thus $n=3,2,1$.
\end{setup}

\begin{setup}\label{n=3}{\bf Three basket case.} First we consider
the case $n=3$. Assume that the basket $Q_i$ is of the type
$\frac{1}{r_i}(a_i,-a_i,1)$ with $a_ib_i\equiv 1$ (mod
$r_i$) and $0<b_i<r_i$ for $i=1,2,3$.  Since $K^3>0$, one
has:
$$1=P_2(X)> \sum_{i=1}^3 \frac{b_i(r_i-b_i)}{2r_i}\geq
\sum_{i=1}^3 \frac{r_i-1}{2r_i}
$$
and so
$$\frac{1}{r_1}+\frac{1}{r_2}+\frac{1}{r_3} >1 .\eqno{(4.2)}$$
One may assume $r_1\leq r_2\leq r_3$. Then
clearly, the only possible solution  for $(r_1,r_2,r_3)$  are
$(2,3,3),(2,3,4),(2,3,5)$ and $(2,2,r_3)$.

\medskip

({\ref{n=3}.1.}) The case $(r_1,r_2,r_3)=(2,2,r_3)$.\\
By (4.1), we have $b'_1=b'_2=b'_3=1$ and $K_X^3=\frac{1}{r_3}$.
Hence $b_1=b_2=1$, and $b_3=1$ or $r_3-1$. Easy computation shows
that $P_4(X)=4 $ (resp. $=5$) if $r_3 \ge 3$ (resp. $r_3=2$). And
also $P_5(X)=6$ (resp. $=7,=9$) if $r_3 \ge 4$ (resp. $r_3=3,=2$).

\medskip

({\ref{n=3}.2.}) The case $(r_1,r_2,r_3)=(2,3,3)$.\\
Computation shows that $K_X^3=\frac{1}{6}$ and $P_4=3, P_5=5$.

\medskip

({\ref{n=3}.3.}) The case $(r_1,r_2,r_3)=(2,3,4)$.\\
Since $b'_1=b'_2=b'_3=1$. Then one gets  a possible case:
\begin{quote}
{\bf (C1).} $(r_1,r_2,r_3)=(2,3,4)$, $K^3=\frac{1}{12}$, $P_2=1$,
$P_3=2$, $P_4=3$ and $P_5=4$.
\end{quote}

({\ref{n=3}.4.}) The case $(r_1,r_2,r_3)=(2,3,5)$.\\
Similarly, by $b'_1=b'_2=b'_3=1$.  So we have found another
possible case:
\begin{quote}
{\bf (C2).} $(r_1,r_2,r_3)=(2,3,5)$, $b_3=1$ or $4$,
$K^3=\frac{1}{30}$, $P_2=1$, $P_3=2$, $P_4=3$ and $P_5=4$.
\end{quote}
\medskip
\end{setup}
\medskip

\begin{setup}\label{n=2}{\bf Two basket case.} Consider the
case $n=2$. We may assume $r_1\leq r_2$.  Also recall that
$b'_1+b'_2=3$ by (4.1). We will distinguish the following two
cases.
\medskip

({\ref{n=2}.1.}) $b'_1=1, b'_2=2$.\\
By $(4.1)$, $1+ K_X^3= \frac{1}{r_1}+\frac{4}{r_2}$. Hence we have
$\frac{5}{r_1} \ge \frac{1}{r_1}+\frac{4}{r_2} >1$. It follows
that $r_1 <5$.

If $r_1=2$, then one gets $\frac{4}{r_2}= K^3+\frac{1}{2}$
and hence $r_2<8$. Noting $(b_2,r_2)=1$ and $b'_2 \le
\frac{1}{2}r_2$, one sees that $r_2=5,7$. Whenever
$(r_1,r_2)=(2,5)$, then computation shows that $K_X^3=
\frac{3}{10}, P_4=4, P_5=7$. Whenever $(r_1,r_2)=(2,7)$, we
have found the possible case:
\begin{quote}
{\bf (C3).} $(r_1,r_2)=(2,7)$, $b_2=2$ or $5$, $K^3=\frac{1}{14}$,
$P_2=1$, $P_3=2$, $P_4=3$ and $P_5=4$.
\end{quote}
\medskip

If $r_1=3$, then $ \frac{4}{r_2}=\frac{2}{3}+K_X^3 > \frac{2}{3}$.
This gives $r_2<6$. The only possibility is $(r_1,r_2)=(3,5)$
since $2=b'_2 \le \frac{1}{2}r_2$ and $(b'_2,r_2)=1$. Computation
shows that  $K_X^3=\frac{2}{15}, P_4=3$ and $P_5=5$.

If $r_1=4$, then similarly we have  $r_2\leq 5$. The only
possibility is $r_2=5$. So $K^3=\frac{1}{20}$ and $P_3=2$,
$P_4=3$, $P_5=4$. We have found the possible case:

\begin{quote}
{\bf (C4).} $(r_1,r_2)=(4,5)$, $b_2=2$ or $3$, $K^3=\frac{1}{20}$,
$P_2=1$, $P_3=2$, $P_4=3$ and $P_5=4$.
\end{quote}
\medskip

({\ref{n=2}.2.}) $b'_1=2, b'_2=1$.\\
By $(4.1)$, $1+ K_X^3= \frac{4}{r_1}+\frac{1}{r_2}$. Hence we have
$\frac{5}{r_1} \ge \frac{1}{r_1}+\frac{4}{r_2} >1$. It follows
that $r_1 <5$. However, $2=b'_1 \le \frac{1}{2}{r_1}$ and
$(2,r_1)=1$ gives a contradiction.

\end{setup}

\begin{setup}\label{n=1}{\bf One basket case.} By (4.1), one has
$b'=3$ and $\frac{9}{r}=1+ K_X^3 >1$. Hence $r<9$. Moreover
$b' \le \frac{r}{2}$ and $(b',r)=1$, so it follows that
$r=7,8$.

If $r=7$,  one gets  $K_X^3=\frac{2}{7},P_4=4, P_5=7$.

If  $r=8$, one gets  $K_X^3=\frac{1}{8}, P_4=3, P_5=5$.
\end{setup}

We summarize all the possible cases with $P_g=1, P_5 <5$:

\begin{cor}\label{class} Let $X$ be a minimal projective
3-fold of general type with $\chi(\OO_X)=0$, $P_g(X)=1$ and
$P_5(X)=4$. Then $X$ has at most 3 baskets of singularities of
type $\frac{1}{r}(a,-a,1)$ and one of the following 4 situations
occurs:
\begin{itemize}
\item [{\bf (C1).}] $(r_1,r_2,r_3)=(2,3,4)$, $K^3=\frac{1}{12};$
\item [{\bf (C2).}] $(r_1,r_2,r_3)=(2,3,5)$, $b_3=1$ or $4$,
$K^3=\frac{1}{30};$
\item [{\bf (C3).}] $(r_1,r_2)=(2,7)$, $b_2=2$
or $5$, $K^3=\frac{1}{14};$
\item [{\bf (C4).}]
$(r_1,r_2)=(4,5)$, $b_2=2$ or $3$, $K^3=\frac{1}{20}.$
\end{itemize}
\end{cor}

Example \ref{1/30} shows that the situation {\bf {(C2)}} really does
occur. We give another example to show the existence of {\bf (C3),
(C4)}.

\begin{exmp} (1) (\cite{C-R}, p153)  The canonical hypersurface
$$X_{12,15}\subset {\bP}(1,3,4,5,6,7)$$ has two terminal
singularities: $1\times \frac{1}{7}(4,-4,1)$, $1\times
\frac{1}{2}(1,-1,1)$. The canonical volume is
$\frac{1}{14}$. This example corresponds to {\bf (C3)}.

(2) (\cite{C-R}, p151) The canonical hypersurface
$$X_{21}\subset {\bP}(1,3,4,5,7)$$
has two terminal singularities: $1\times
\frac{1}{4}(1,-1,1)$, $1\times \frac{1}{5}(3,-3,1)$. The
canonical volume is $\frac{1}{20}$. This example
corresponds to {\bf (C4)}.
\end{exmp}

It is interesting to ask:
\begin{question} Does {\bf (C1)} really occur?
\end{question}

\begin{setup}{\bf Proof of Theorem \ref{main}}.
\begin{proof} Let $X$ be a minimal projective 3-fold of general
type (admitting at worst $\bQ$-factorial terminal singularities)
with $\chi(\OO_X)\leq 0$. Recall that one has
$$\chi(\OO_X)=1-q+h^2(\OO_X)-p_g$$
where the irregularity $q:=h^1(\OO_X)$ and the geometric
genus $p_g:=h^3(\OO_X)$. Since $\Vol(X)\geq \frac{1}{3}$
whenever $p_g\geq 2$ by \cite{mathann}, and $\Vol(X)\geq
\frac{1}{22}$ whenever $q>0$ by Corollary \ref{q>0}, we may
assume, from now on, that $p_g\leq 1$ and $q=0$. Therefore
the assumption $\chi(\OO_X)\leq 0$ implies $p_g=1$,
$h^2(\OO_X)=0$ and finally $\chi(\OO_X)=0$.

Whenever $P_5(X)\geq 5$, Theorem \ref{p3} (ii) and (iii)
says $\Vol(X)>\frac{1}{25}$.

Whenever $P_5(X)\leq 4$, $p_g>0$ and $\chi(\OO_X)=0$,
Corollary \ref{class} says that $\Vol(X)\geq \frac{1}{30}$.
Furthermore $\Vol(X)=\frac{1}{30}$ implies that $X$
corresponds exactly to the situation ${\bf (C2)}$ in the
list of Corollary \ref{class}. This completes the proof.
\end{proof}
\end{setup}



\begin{thebibliography}{99}
\bibitem{At} M. Atiyah, {\em Vector bundles over an elliptic
curves}, Proc. London Math. Soc. {\bf 7} (1957), 414-452.

\bibitem{Bom} E. Bombieri, {\em
Canonical models of surfaces of general type}. Inst. Hautes Etudes
Sci. Publ. Math. {\bf 42} (1973), 171--219.


\bibitem{JC-H-mz} J. A. Chen, C. D. Hacon, {\em Pluricanonical
systems on irregular 3-folds of general type}. Math. Z. {\bf
255} (2007), no. 2, 343-355

\bibitem{mathann} M. Chen, {\em A sharp lower bound for the
canonical volume of 3-folds of general type}, Math. Ann. {\bf
337} (2007), 165-181.

\bibitem{JMSJ3} M. Chen, {\em Inequalities of Noether type
for 3-folds of general type}, J. Math. Soc. Japan {\bf 56} (2004),
1131-1155.

\bibitem{MPCPS} M. Chen, {\em Canonical stability in terms of
singularity index for algebraic threefolds}, Math. Proc. Camb.
Phil. Soc. {\bf 131} (2001), 241-264.

\bibitem{Chen-Zuo} M. Chen, K. Zuo, {\em Complex projective
threefolds with non-negative canonical Euler-Poincare
characteristic}, Communications in Analysis and Geometry {\bf 16}
(2008), 159-182.

\bibitem{Ci} C. Ciliberto, {\em The bicanonical map for surfaces of general type},
Proc. Symposia in Pure Math. {\bf 62} (1997), 57-83.

\bibitem{C-R} A. Corti, M. Reid, {\em Explicit birational geometry of 3-folds}.
London Mathematical Society, Lecture Note Series, 281. Cambridge
University Press, Cambridge, 2000.

\bibitem{H-M} C. D. Hacon and J. M$^{\text{\rm c}}$Kernan, {\em
Boundedness of pluricanonical maps of varieties of general type},
Invent. Math. {\bf 166} (2006), 1-25.

\bibitem{Fletcher} A. R. Iano-Fletcher,
{\em Inverting Reid's exact plurigenera formula}. Math. Ann. {\bf 284}
(1989), no. 4, 617-629.
\bibitem{KV} Y. Kawamata,
{\em A generalization of Kodaira-Ramanujam's vanishing theorem},
Math. Ann. {\bf 261} (1982), 43-46.

\bibitem{KMM} Y. Kawamata, K. Matsuda, K. Matsuki, {\em Introduction to the
minimal model problem}, Adv. Stud. Pure Math. {\bf 10} (1987),
283-360.
\bibitem{K-M} J. Koll\'ar, S. Mori, Birational geometry of algebraic
varieties, 1998, Cambridge Univ. Press.

\bibitem{ohno} K. Ohno, {\em Some inequalities for minimal
fibrations of surfaces of general type over curves}, J.
Math. Soc. Japan {\bf 44} (1992), 643-666.


\bibitem{YPG} M. Reid,  {\em Young person's guide to canonical
singularities}, Proc. Symposia in pure Math. {\bf 46} (1987),
345-414.
\bibitem{Reid83} M. Reid, {\em Minimal models of
canonical 3-folds}, Adv. Stud. Pure Math. {\bf 1} (1983), 131-180.
\bibitem{Reider} I. Reider, {\em Vector bundles of rank 2 and linear systems on
algebraic surfaces}, Ann. Math. {\bf 127} (1988), 309-316.

\bibitem{Tak} S. Takayama, {\em Pluricanonical systems on algebraic
varieties of general type}, Invent. Math. {\bf 165} (2006), 551
-- 587.


\bibitem{Tsuji} H. Tsuji, {\em
Pluricanonical systems of projective varieties of general type, I}, Osaka J. Math. {\bf 43}
(2006), 967--995.

\bibitem{V} E. Viehweg, {\em Vanishing theorems},
J. reine angew. Math. {\bf 335} (1982), 1-8.

\bibitem{V1}
E. Viehweg, {\em Weak positivity and the additivity of the
Kodaira dimension for certain fibre spaces}. Proc.
Algebraic Varieties and Analytic Varieties, Tokyo 1981.
Adv. Studies in Math. {\bf 1}, Kinokunya-North-Holland
Publ. 1983, 329-353

\end{thebibliography}
\end{document}